\documentclass[12pt,a4paper]{amsart}
\usepackage[all]{xy}
\usepackage{amscd}
\usepackage{amsfonts}
\usepackage{amsmath}
\usepackage{amssymb}
\usepackage{amsthm}
\usepackage[colorlinks,linkcolor=blue]{hyperref}
\usepackage{calc}
\usepackage{extarrows}
\usepackage{graphicx}
\usepackage{graphics}
\usepackage[margin=2.9cm]{geometry}
\usepackage{mathrsfs}

\theoremstyle{plain}
\newtheorem{thm}{Theorem}[subsection]
\newtheorem{lema}[thm]{Lemma}
\newtheorem{cor}[thm]{Corollary}
\newtheorem{ppst}[thm]{Proposition}

\theoremstyle{definition}
\newtheorem{defn}[thm]{Definition}
\newtheorem{rem}[thm]{Remark}
\newtheorem{exam}[thm]{Example}

\newcommand{\bigboxplus}{
  \mathop{
    \vphantom{\bigoplus}
    \mathchoice
      {\vcenter{\hbox{\resizebox{\widthof{$\displaystyle\bigoplus$}}{!}{$\boxplus$}}}}
      {\vcenter{\hbox{\resizebox{\widthof{$\bigoplus$}}{!}{$\boxplus$}}}}
      {\vcenter{\hbox{\resizebox{\widthof{$\scriptstyle\oplus$}}{!}{$\boxplus$}}}}
      {\vcenter{\hbox{\resizebox{\widthof{$\scriptscriptstyle\oplus$}}{!}{$\boxplus$}}}}
  }\displaylimits
}

\newcommand{\sign}{\textup{sign}}

\title[]{Examples on the sharpness of an inequality about multiplicities over hyperfields}
\author{Ziqi Liu}
\address{Jilin University, Changchun, Jilin, China}
\email{liuzq0616@mails.jlu.edu.cn}

\thanks{This work is an undergraduate research under the instruction of Professor Matthew Baker in Georgia Institute of Technology. The author also thanks professor Oliver Lorscheid for helpful discussions and feedback in October, 2019. The main part of this work was finished in December, 2019 and the author made some modifies recently.
}

\begin{document}

\begin{abstract}
  In this paper, the author introduces hyperfields and give some facts about the roots and multiplicities of polynomials over hyperfield based on \cite{2} and \cite{3}. Then he tests the sharpness of an inequality in Baker's former work \cite{2} about the behavior of multiplicities under homomorphisms between hyperfields and show that the inequality is not sharp under the natural homomorphisms $\mathbb{C}\rightarrow\mathbb{P}$ and $\mathbb{C}\rightarrow\mathbb{V}$ while it is sharp under the natural homomorphisms $\mathbb{C}\rightarrow\mathbb{K}$, $\mathbb{R}\rightarrow\mathbb{S}$ and $\mathbb{R}\rightarrow\mathbb{T}$ according to the previous work \cite{2} of Baker and Lorscheid.
\end{abstract}

\maketitle
\tableofcontents

\newpage
\section{Some Backgrounds}
\subsection{The Definition of Hyperfields}
\begin{defn}
A \textbf{hyperoperation} on a set $S$ is a map $\boxplus:S\times S\rightarrow 2^S\backslash\{\varnothing\}$. For given hyperoperation $\boxplus$ on $S$ and non-empty subsets $A,B$ of $S$, define
$$A\boxplus B:=\bigcup_{a\in A,b\in B}(a\boxplus b)$$
A hyperoperation $\boxplus$ is \textbf{commutative} if $a\boxplus b=b\boxplus a$ for all $a,b\in S$. If not especially mentioned, hyperoperations in this paper will always be commutative.\\
A hyperoperation $\boxplus$ is \textbf{associative} if $a\boxplus(b\boxplus c)=(a\boxplus b)\boxplus c$ for all $a,b,c\in S$.
\end{defn}
\begin{defn}
Given an associative hyperoperation $\boxplus$, the \textbf{hypersum} of $x_1,x_2,\dots,x_m$ for $m\geq2$ is recursively defined as
$$x_1\boxplus\cdots\boxplus x_n:=\bigcup_{x'\in x_2\boxplus\cdots\boxplus x_n}x_1\boxplus x'$$
\end{defn}
\begin{defn}
A \textbf{hypergroup} is a tuple $(G,\boxplus,0)$, where $\boxplus$ is an associative hyperoperation on $G$ such that:\\
(1) $0\boxplus x=\{x\}$ for all $x\in G$;\\
(2) For every $x\in G$ there is a unique element $-x$ of $G$ such that $0\in x\boxplus-x$;\\
(3) $x\in y\boxplus z$ if and only if $z\in x\boxplus(-y)$.\\
Here $-x$ is often called as the \textbf{hyperinverse} of $x$ and (3) as the reversibility axiom.
\end{defn}
\begin{defn}
A (Krasner) \textbf{hyperring} is a tuple $(R,\odot,\boxplus,1,0)$ such that:\\
(1) $(R,\odot,1)$ is a commutative monoid;\\
(2) $(R,\boxplus,0)$ is a commutative hypergroup;\\
(3) $0\odot x=x\odot0=0$ for all $x\in R$;\\
(4) $a\odot(x\boxplus y)=(a\odot x)\boxplus(a\odot y)$ for all $a,x,y\in R$;\\
(5) $(x\boxplus y)\odot a=(x\odot a)\boxplus(y\odot a)$ for all $a,x,y\in R$.\\
In the following part, we often use the underlying set $R$ to refer to a hyperring and may omit $\odot$ if there is no likehood of confusion.
\end{defn}
\begin{defn}
A hyperring $F$ is called a \textbf{hyperfield} if $0\neq 1$ and every non-zero element of $F$ has a multiplicative inverse.
\end{defn}

\begin{exam}
If $(\mathbb{F},\cdot,+)$ is a field, then $\mathbb{F}$ can be trivially associated with a hyperfield $(\mathbb{F},\odot,\boxplus)$ where $x\odot y=x\cdot y$ and $x\boxplus y=\{x+y\}$ for all $x,y\in\mathbb{F}$.\\
In the following context, when we mention a field $\mathbb{F}$, we may actually refer to the hyperfield associated with $\mathbb{F}$.
\end{exam}
\begin{exam}
Consider $\mathbb{K}=(\{0,1\},\odot,\boxplus,1,0)$ with the usual multiplication rule and a hyperaddition $\boxplus$ defined by
$$0\boxplus0=\{0\},\qquad1\boxplus0=0\boxplus1=\{1\},\qquad 1\boxplus1=\{0,1\}$$
then $\mathbb{K}$ is a hyperfield, called the \textbf{Krasner hyperfield}.
\end{exam}
\begin{exam}
Consider $\mathbb{S}=(\{0,1,-1\},\odot,\boxplus,1,0)$ with the usual multiplication rule and a hyperaddition $\boxplus$ generated by
$$x\boxplus x=\{x\},\quad x\boxplus 0=\{x\},\quad1\boxplus-1=\{-1,0,1\}$$
then $\mathbb{S}$ is a hyperfield, called the \textbf{hyperfield of signs}.
\end{exam}
\begin{exam}
Let $\mathbb{T}:=\mathbb{R}\cup\{-\infty\}$ as sets and define hyperoperation $\boxplus$ as
$$ x\boxplus y=\left\{
\begin{aligned}
\{\max\{x,y\}\} &, & x\neq y \\
\{z\in\mathbb{T}:z\leq x\} &,  & x=y
\end{aligned}
\right.
$$
and $\odot$ as $x\odot y=x+y$. Then $\mathbb{T}$ is a hyperfield, called the \textbf{tropical hyperfield}.
\end{exam}
\begin{rem} More generally, let $\Gamma$ be a totally ordered abelian group (written multiplicatively) one can define a canonical hyperfield structure on set $\Gamma\cup\{0\}$ where
\begin{itemize}
  \item the multiplication $\odot$ is multiplication in $\Gamma$ with $0\odot x=0$ for all $x\in\Gamma\cup\{0\}$
  \item the hyperaddition $\boxplus$ is defined as $x\boxplus x:=\{y:y\leq x\}$ and $x\boxplus y:=\max\{x,y\}$ for $x\neq y$.
\end{itemize}
with $0\leq x$ for all $x\in\Gamma$. One such hyperfield is a \textbf{valuative hyperfield}.
\end{rem}
\begin{exam}
Let $\mathbb{P}=S^1\cup\{0\}$, where $S^1=\{z\in\mathbb{C}:|z|=1\}$ is the complex unit circle. Then one can define a hyperfield structure on $\mathbb{P}$ where the multiplication is the usual one and the hyperaddition is defined as
$$ x\boxplus y=\left\{
\begin{array}{cll}
\{x\}, & & y=0 \\
\{0,x,-x\}, & & y=-x \\
\{\frac{\alpha x+\beta y}{||\alpha x+\beta y||}:\alpha,\beta\in\mathbb{R}_+\}, &  & \textup{otherwise}
\end{array} \right. $$
This hyperfield structure on $S^1\cup\{0\}$ is called the \textbf{phase hyperfield}.
\end{exam}
\begin{exam}
Let $\mathbb{V}$ be the set $\mathbb{R}_{\geq0}=\mathbb{R}_+\cup\{0\}$ of nonnegative real numbers with the usual multiplication and the hyperaddition is defined as
$$x\boxplus y=\{z\in\mathbb{R}_{\geq0}:|x-y|\leq z\leq x+y\}$$
Then $\mathbb{V}$ is a hyperfield named the \textbf{Viro hyperfield} (or the triangle hyperfield).
\end{exam}
\begin{defn}
Let $\mathbb{F}_1,\mathbb{F}_2$ be two hyperfields. A map $f:\mathbb{F}_1\rightarrow\mathbb{F}_2$ is called a \textbf{hyperfield homomorphism} if for all $x,y\in \mathbb{F}_1$\\
(1) $f(0)=0,f(1)=1$;\\
(2) $f(xy)=f(x)f(y)$ and $f(x\boxplus y)\subseteq f(x)\boxplus f(y)$.
\end{defn}
\begin{exam}\label{RS}
Consider the field of real numbers $\mathbb{R}$ and the hyperfield of signs $\mathbb{S}$, a natural hyperfield homomorphism can be defined as
\begin{center}
$f:\mathbb{R}\longrightarrow\mathbb{S},\quad x\longmapsto \sign(x)$
\end{center}
where $\sign:\mathbb{R}\rightarrow\{1,-1,0\}$ is the function which maps a real number to its sign.
\end{exam}
\begin{exam}\label{CV}
Consider the field of complex numbers $\mathbb{C}$ and the Viro hyperfield $\mathbb{V}$, one can see a natural hyperfield homomorphism
\begin{center}
$f:\mathbb{C}\longrightarrow\mathbb{V},\quad re^{i\theta}\longmapsto r$
\end{center}
where $r\geq0$ and $\theta\in[0,2\pi]$.
\end{exam}
\begin{ppst}\label{lemma}
Given two hyperfields $\mathbb{F}_1$ and $\mathbb{F}_2$ and a homomorphism $f:\mathbb{F}_1\rightarrow\mathbb{F}_2$, then for any $x_1,x_2,\dots,x_n$ in $\mathbb{F}_1$, we have
\begin{center}
$f(x_1\boxplus x_2\boxplus\cdots\boxplus x_n)\subseteq f(x_1)\boxplus f(x_2)\boxplus\cdots\boxplus f(x_n)$
\end{center}
holds, where $n\geq2$ is an arbitrary integer.

\begin{proof}
We prove it by induction. The $n=2$ case is obviously valid by definition.\\
Assume that our claim is valid for $n-1$, then for case $n$, we have
$$\begin{aligned}
f(x_1\boxplus x_2\boxplus\cdots\boxplus x_n)&=f(\bigcup_{x\in x_2\boxplus\cdots\boxplus x_n} x_1\boxplus x)\\
&=\bigcup_{x\in x_2\boxplus\cdots\boxplus x_n}f(x_1\boxplus x)\\
&\subseteq\bigcup_{x\in x_2\boxplus\cdots\boxplus x_n}f(x_1)\boxplus f(x)\\
&\subseteq\bigcup_{y\in f(x_2\boxplus\cdots\boxplus x_n)}f(x_1)\boxplus y\\
&\subseteq\bigcup_{y\in f(x_2)\boxplus\cdots\boxplus f(x_n)}f(x_1)\boxplus y\\
&=f(x_1)\boxplus f(x_2)\boxplus\cdots\boxplus f(x_n)
\end{aligned}$$
where we use our inductive assumption for the last $\subseteq$.\\
Hence our claim is true for $n$ and then we finish the induction and the whole proof.
\end{proof}
\end{ppst}

\subsection{The Definition of Multiplicities over Hyperfields}
\begin{defn}
Given a hyperfield $\mathbb{F}$, a \textbf{polynomial} over $\mathbb{F}$ (or with coefficients in $\mathbb{F}$) is a map $p:\mathbb{F}\rightarrow2^{\mathbb{F}}$ that
\begin{center}
$a\longmapsto c_na^n\boxplus c_{n-1}a^{n-1}\boxplus\cdots\boxplus c_1a\boxplus c_0$
\end{center}
where $\{c_i\}_{i=0}^n\subset\mathbb{F}$ and $c_n\neq0$.\\
For such $p$, denote it by $p(T)=c_nT^n+c_{n-1}T^{n-1}+\cdots+c_1T+c_0$.
\end{defn}
\begin{exam}
The polynomial $p(T)=1T^3-2$ over the tropical hyperfield $\mathbb{T}$ is
$$p(a)=1a^3\boxplus(-2)=\left\{\begin{aligned}
-2\quad&,\,\,a<-1\\
[-\infty,-2]&,\,\,a=-1\\
1a^3\quad&,\,\,a>-1
\end{aligned}\right.$$
where $<$ is the natural order in the real axis and $1a^3$ means $1\odot a\odot a\odot a$ in $\mathbb{T}$.
\end{exam}
\begin{ppst}
Let $\mathbb{F}$ be a hyperfield, the set of all polynomials over $\mathbb{F}$ can be naturally endowed with two hyperoperations\\
(1) $p\boxdot q=e_{mn}T^{mn}+\cdots+e_1T_1+e_0$ where $e_i=\mathop{\bigboxplus_{\mathbb{F}}}\limits_{k+l=i}c_kd_l$;\\
(2) $p\boxplus q=e_{n}T^{n}+\cdots+e_1T_1+e_0$ where $e_i=c_i\boxplus_{\mathbb{F}} d_i$ for $i\leq m$ and $e_i=c_i$ for $i>m$.\\
for $p(T)=c_nT^n+\cdots+c_1T+c_0$ and $q(T)=d_mT^m+\cdots+d_1T+d_0$ with $n\geq m$.\\
Therefore, one can get a new hyperstructure over that set, denoted by $\textup{Poly}(\mathbb{F})$.
\end{ppst}
\begin{rem} In \cite{3}, this hyperstructure is called a polynomial hyperring (but it is clearly not a hyperring). In some materials like \cite{1}, this hyperstructure is identified as a new concept: superring.
\end{rem}
\begin{defn}
Let $p(T)=c_n T^n+ c_{n-1}T^{n-1}+\cdots+ c_1T+ c_0$ be a polynomial over a hyperfield $\mathbb{F}$, an element $a\in\mathbb{F}$ is called a \textbf{root} of $p$ if and only if either the following equivalent conditions is satisfied:\\
(1) $0\in p(a)=c_n a^n\boxplus c_{n-1}a^{n-1}\boxplus \cdots\boxplus c_1a\boxplus c_0$;\\
(2) there exist $q$ in Poly$(\mathbb{F})$ such that $p\in(T-a)\boxdot q$ (the symbol $\boxdot$ is sometimes omitted).
\end{defn}
\begin{defn}
Let $\displaystyle p(T)=c_n T^n+ c_{n-1}T^{n-1}+\cdots+ c_1T+ c_0$ be a polynomial over a hyperfield $\mathbb{F}$, if $a$ is not a root of $p$, set mult$_{a}(p)=0$. If $a$ is a root of $p$, define
\begin{center}
mult$_{a}(p)$ = $1+\max\{\textup{mult}_{a}(q):p\in (T-a)q\}$
\end{center}
as the \textbf{multiplicity} of the root $a$ of $p$. Moreover, for a nonempty set $S$, define
\begin{center}
mult$_{S}(p)$ = $1+\max\{\textup{mult}_{S}(q):p\in (T-a)q$ for some $a\in S\}$
\end{center}
It is clearly that mult$_S(p)\leq\textup{deg}(p)$ for any $S\subset\mathbb{F}$.
\end{defn}
\begin{exam} Given a polynomial $p(T)=T^2+3T+1$ over the Viro hyperfield $\mathbb{V}$.\\
(1) Since $p(3)=(3^2)\boxplus(3^2)\boxplus1=9\boxplus9\boxplus1=[0,19]$
one has $0\in p(3)$ and hence mult$_3(p)\geq1$. From
$$(T+3)\boxdot(T+3)=\{T^2+aT+9:a\in[0,6]\}$$
one can see $p(T)\notin(T+3)(T+3)$ and then mult$_3(p)=1$.\\
(2) Given $S=[1,+\infty)$, from (1) one can see mult$_S(p)\geq1$. Suppose mult$_S(p)=2$, then there exist an $s\in S$ such that
$$p(T)\in(T+s)\boxdot(T+s)=\{T^2+aT+s^2:a\in[0,2s]\}$$
In this case, $s^2=1$ and thus $s=1$, which is impossible. Therefore, mult$_S(p)=1$.
\end{exam}
\begin{defn}
Given hyperfields $\mathbb{F}_1,\mathbb{F}_2$ and a polynomial $\overline{p}$ over $\mathbb{F}_2$ if there exists a hyperfield homomorphism $f:\mathbb{F}_1\rightarrow\mathbb{F}_2$, then a \textbf{corresponding polynomial} of $\overline{p}$ with respect to $f$ is a polynomial $p$ in set
$$\{p(T)=c_{n}T+\cdots+c_1T+c_0\,|\,f(c_{n})T+\cdots+f(c_1)T+f(c_0)=\overline{p}(T)\}$$
which is a subset of Poly$(\mathbb{F}_1)$.\\
Also, for $p$ in Poly$(\mathbb{F}_1)$, one also denotes $f(c_{n})T+\cdots+f(c_1)T+f(c_0)$ by $\overline{p}$.
\end{defn}
\begin{exam} It is natural for us to consider if multiplicity of a polynomial can be preserved by hyperfield homomorphisms. Following are two simple examples.\\
(1) Consider $p(T)=T^2-3T+2$ over the field of real numbers $\mathbb{R}$ and the hyperfield homomorphism $f:\mathbb{R}\rightarrow\mathbb{S}$ described in Example \ref{RS}, one is able to see that the corresponding polynomial $\overline{p}(T)=T^2-T+1$ with respect to $f$ satisfies
\begin{center}
mult$\displaystyle_1(\overline{p})=2=\sum_{a>0}\textup{mult}_a(p)=\sum_{a\in f^{-1}(1)}\textup{mult}_a(p)$.
\end{center}
(2) Consider $q(T)=T^3+2T-3$ over the complex field $\mathbb{C}$ and the hyperfield homomorphism $g:\mathbb{C}\rightarrow\mathbb{R}$ described in Example \ref{CV}, it is not hard to check that the corresponding  polynomial $\overline{q}(T)=T^3+4T+5$ with respect to $g$ satisfies
\begin{center}
mult$\displaystyle_2(\overline{q})=2>0=\sum_{|a|=2}\textup{mult}_a(q)=\sum_{a\in g^{-1}(2)}\textup{mult}_a(q)$.
\end{center}
For (2), we have mult$_2(\overline{q})\leq2$ and
\begin{center}
$T^3+4T+5\in(T+2)(T^2+2T+2.5),\quad T^2+2T+2.5\in(T+2)(T+1.25)$
\end{center}
where one should notice that $-x=x$ holds for each $x\in\mathbb{V}$.
\end{exam}
\begin{ppst}\label{power}
For a monic polynomial $p(T)=T^n+c_{n-1}T^{n-1}+\cdots+c_0$ over a hyperfield $\mathbb{F}$ with $x\boxplus x=\{x\}$, the following two statements are equivalent:\\
(1) $a$ is a root of $p$ with mult$_a(p)=n$;\\
(2) $c_{n-k}=(-a)^k$ for $k=1,\dots,n$.

\begin{proof}
We prove our claim by induction on $n$. First, it is clear that our claim is true for $n=1$ and we know it is true for $n=2$ since
$$(T-a)\boxdot(T-a)=\{T^2+(-a)\boxplus(-a)T+(-a)^2\}=\{T^2+(-a)T+(-a)^2\}$$
Suppose our claim is valid for $n$, consider the case $n+1$.\\
(1)$\Rightarrow$(2). If $a$ is a root of $p(T)=T^{n+1}+c_nT^n+\cdots+c_0$ with mult$_a(p)=n+1$, then there exist a polynomial $q(T)$ such that mult$_a(q)=n$ and $p\in(T-a)q$. By our inductive assumption, one knows that
$$q(T)=T^n-aT^{n-1}+\cdots+(-a)^{n-1}T+(-a)^n$$
Therefore, from $p\in(T-a)q$ one has
\begin{center}
$c_0=(-a)^{n+1}$ and $c_{n+1-k}\in (-a)(-a)^{k-1}\boxplus(-a)^k={(-a)^k}$
\end{center}
this part is done.\\
(2)$\Rightarrow$(1). If $c_{n+1-k}=(-a)^k$ for $k=1,\dots,n+1$, one can get $p\in(T-a)q$ where
$$q(T)=T^n+(-a)T^{n-1}+\cdots+(-a)^{n-1}T+(-a)^n$$
By our inductive assumption, we have mult$_a(q)=n$ and it follows immediately that mult$_a(p)=n+1$.\\
In conclusion, our claim is still true for $n+1$ and then we are done here.
\end{proof}
\end{ppst}

\section{The Inequality about Multiplicities over Hyperfields}
\subsection{Introduction to the Inequality}
\begin{ppst}
Given two hyperfields $\mathbb{F}_1,\mathbb{F}_2$ and a homomorphism $f:\mathbb{F}_1\rightarrow\mathbb{F}_2$. Let
\begin{center}
$p(T)=T^n+ c_{n-1}T^{n-1}+\cdots+c_1T +c_0$
\end{center}
be a polynomial with coefficients in a hyperfield $\mathbb{F}_1$, if $\alpha$ is a root of $p$, then $f(\alpha)$ is a root of
\begin{center}
$\overline{p}(T)=T^n+ f(c_{n-1})T^{n-1}+\cdots+f(c_1)T+ f(c_0)$
\end{center}
with coefficients in $\mathbb{F}_2$.
\end{ppst}
This proposition is a corollary of Proposition \ref{lemma} and one can then obtain the following theorem.
\begin{thm}\label{tar}
Given a field $\mathbb{L}$, a hyperfield $\mathbb{F}$ and a homomorphism $f:\mathbb{L}\rightarrow\mathbb{F}$, we have
$$\textup{mult}_b(\overline{p})\geq\sum_{a\in f^{-1}(b)}\textup{mult}_a(p),\quad\textup{for all}\,\,b\in\mathbb{F}$$
for every $\overline{p}$ and its corresponding polynomial $p$ with respect to $f$.\\
Moreover, if $\sum_{b\in\mathbb{F}}\textup{mult}_b(\overline{p})\leq\textup{deg}(\overline{p})$ and $p$ splits into a product of linear factors over $\mathbb{L}$, then the equality holds in this theorem.
\end{thm}
This theorem is Proposition B in \cite{3}. In fact, we can write a more general theorem with homomorphisms from hyperfield $\mathbb{L}$ to $\mathbb{F}$ where polynomials in Poly$(\mathbb{L})$ have unique decomposition. In the following part, We will focus on whether the inequality in Theorem \ref{tar} is sharp under different hyperfield homomorphisms, that is, whether we can always find a $p$ for each $\overline{p}$ such that the inequality holds.

\newpage
\begin{ppst}\label{2.1.3}
Let $p(T)=T^n+ c_{n-1}T^{n-1}+\cdots+ c_1T+ c_0$ be a polynomial with coefficients in a hyperfield $\mathbb{F}$, if mult$_a(p)=n$, then we have $c_n=(-a)^n$.

\begin{proof}
We will prove it by induction, since our definition of multiplicity is inductive.\\
For $n=1$, if $0\in a\boxplus c_0$, then $c_0\in 0\boxplus(-a)=\{-a\}$ which implies that $c_0=-a$.\\
Assume that our claim is true for $n-1$, we want to prove that it is also true for $n$.\\
Let $a$ be a root of polynomial $p(T)$ with mult$_a(p)=n$, then from
\begin{center}
mult$_{a}(p)$ = $1+\max\{\textup{mult}_{a}(q):p\in (T-a)q\}$
\end{center}
we know that $\max\{\textup{mult}_{a}(q):p\in (T-a)q\}=n-1$. So there are some
$$q(T)=T^{n-1}+ d_{n-2}T^{n-2}+\cdots+ d_1T+ d_0$$
such that $a$ is a root of $q(T)$ with mult$_a(p)=n-1$. By our assumption, $d_0=(-a)^n$, hence from $p\in(T-a)q$ we get $c_0=-ad_0=(-a)^n$, our induction is done.
\end{proof}
\end{ppst}

\begin{ppst}\label{2.1.4}
Let $p(T)=c_nT^n+ c_{n-1}T^{n-1}+\cdots+ c_1T+c_0$ be a polynomial over $\mathbb{F}$, then $0$ is a root of $p$ if and only if $c_0=0$.

\begin{proof}
If $0$ is a root of $p$, then there exist $\{d_i\}_{i=0}^{n-1}=\{c_{i+1}\}_{i=0}^{n-1}$ such that
\begin{center}
$c_0=0,\,\,c_i\in 0\boxplus d_{i-1}$ for $i=1,\dots,n-1$ and $c_n=d_{n-1}$
\end{center}
If $c_0=0$, then it is clear that $0\in p(0)=0\boxplus 0\boxplus\cdots\boxplus 0=\{0\}$.
\end{proof}
\end{ppst}

\begin{ppst}\label{2.1.5}
Let $p(T)=T^n+ c_{n-1}T^{n-1}+\cdots+ c_1T+c_0$ be a polynomial with coefficients in $\mathbb{F}$, if $c_0=0$ then the set of nonzero roots of $p$ is the same as the set of nonzero roots of $q$ where $q(T)=T^{n-1}+ c_{n-1}T^{n-2}+\cdots+ c_1$.

\begin{proof}
It is sufficient for us to prove that
$$0\in a^n\boxplus c_{n-1}a^{n-1}\boxplus\cdots\boxplus c_1a\Longleftrightarrow 0\in a^{n-1}\boxplus c_{n-1}a^{n-2}\boxplus\cdots\boxplus c_1$$
for none zero $a\in\mathbb{F}$, which is obvious.\\
Furthermore, the multiplicities of non-zero roots don't change.
\end{proof}
\end{ppst}

\begin{ppst}\label{2.1.6}
Let $p(T)=c_nT^n+ c_{n-1}T^{n-1}+\cdots+c_kT^k$ be a polynomial over $\mathbb{F}$ and for some $0<k<n$, if $c_k\neq0$, then $p$ has 0 as its root with multiplicity of $k$.

\begin{proof}
It is clear that $0$ is a root of $p(T)$ and then we will prove mult$_{0}(p)=k$ by induction.\\
If $k=1$, consider the $ p_1(T)=d_{n-1}T^{n-1}+ d_{n-2}T^{n-2}+\cdots+ d_0$ in
\begin{center}
mult$_{0}(p)$ = $1+\max{\{\textup{mult}_{0}(p_{1}):p\in Tp_{1}\}}$
\end{center}
We have $c_n=d_{n-1}$ and
\begin{center}
$c_i\in 0\boxplus d_{i-1}=\{d_{i-1}\}$ for $i=1,\dots,n-1$
\end{center}
which implies that $p_1(T)\in\{c_nT^{n-1}+ c_{n-1}T^{n-2}+\cdots+ c_1\}$. According to Proposition \ref{2.1.4}, $0$ is not a root of $p_1$ and then mult$_{0}(p)$=1, so our claim is true for $k=1$.\\
Suppose our claim is true for $k-1$, we consider the case with $k$. From our assumption,
\begin{center}
mult$_{0}(p)$ = $k+\max\{\textup{mult}_{0}(p_{k}):p_{k-1}\in Tp_{k}\}$.
\end{center}
where $p_{k-1}(T)=c_nT^{n-k+1}+\cdots+ c_{k}T$. Let $p_{k}(T)=d_{n-k}T^{n-k}+\cdots+ d_{0}$, we have
\begin{center}
$c_{n}=d_{n-k}$ and $c_i\in 0\boxplus d_{i-k}=\{d_{i-k}\}$ for $i=k,\dots,n-1$
\end{center}
which implies that $p_k(T)\in\{c_nT^{n-k}+ c_{n-1}T^{n-k-1}+\cdots+ c_k\}$. According to Proposition \ref{2.1.5}, $0$ is not a root of $p_k$ and then mult$_{0}(p)=k$, so we are done.
\end{proof}
\end{ppst}
\begin{rem}\label{2.1.7}
Propositions \ref{2.1.4}, \ref{2.1.5} and \ref{2.1.6} show that polynomials over hyperfields share some similar nature with respect to the root $0$. Moreover, Proposition \ref{2.1.4} shows that the proposition can be an equivalent statement.\\
Notice that for given hyperfield homomorphism $f:\mathbb{F}_1\rightarrow\mathbb{F}_2$, $f(a)=0$ for some $a\neq0$ implies that $f(ax)\subseteq\{0\}$ for all $x\in\mathbb{F}_1$. Set $x=a^{-1}y$, then $f(y)\subseteq\{0\}$ for all $y\in\mathbb{F}_1$, which implies that $f$ is trivial. Therefore, a nontrivial hyperfield homomorphism preserves $0$ as well as the behaviours of root 0 of polynomials.
\end{rem}
\begin{cor}\label{2.1.8}
For a nontrivial hyperfield homomorphism $f:\mathbb{F}_1\rightarrow\mathbb{F}_2$, if mult$_0(p)=k$ for polynomial $p(T)=T^n+c_{n-1}T^{n-1}+\cdots+c_1T_1+c_0$ over $\mathbb{F}_1$, then we have mult$_0(\overline{p})=k$ for
$$\overline{p}(T)=T^n+f(c_{n-1})+\cdots+f(c_1)T+f(c_0)$$
in \textup{Poly}$(\mathbb{F}_2)$.
\end{cor}
\begin{rem}\label{2.1.9}
This corollary shows that when studying the behaviors of multiplicities of polynomials under nontrivial homomorphisms, one only needs to consider the nonzero roots.
\end{rem}

\subsection{The Sharpness under Natural Hyperfield Homomorphism $\mathbb{C}\rightarrow\mathbb{V}$}
\begin{ppst}\label{2.2.1}
For $x,y\in\mathbb{V}$, $0\in x\boxplus y$ if and only if $x=y$.

\begin{proof}
If $0\in x\boxplus y=\{z:|x-y|\leq z\leq x+y\}$, then $|x-y|=0$ and thus $x=y$.\\
If $x=y$, then $0\in x\boxplus y=x\boxplus x=\{z:0\leq z\leq 2x\}$.
\end{proof}
\end{ppst}

\begin{ppst}\label{2.2.2}
Consider the polynomial $p(T)=T^2+ c_1T+ c_0$ which coefficients in $\mathbb{V}$, then $p$ has a root $a$ with mult$_a(p)=2$ if and only if $c_0=a^2$ and $c_1^2\leq4c_0$.

\begin{proof}
From Proposition \ref{2.1.4}, if $a$ is a root of $p$ with mult$_a(p)=2$, we have $a=\sqrt{c_0}$ and
$$0\in p(\sqrt{c_0})=c_0\boxplus c_1\sqrt{c_0}\boxplus c_0\Longrightarrow c_1\sqrt{c_0}\in c_0\boxplus c_0$$
which gives that $c_1^2\leq 4c_0$.\\
If $c_0=a^2$ and $c_1^2\leq 4c_0$, then we have $p\in(T+a)q$ where $q(T)=T+a$ since
$$c_1\in a\boxplus a=\sqrt{c_0}\boxplus\sqrt{c_0}=[0,2\sqrt{c_0}]$$
Therefore, we are done.
\end{proof}
\end{ppst}
\begin{ppst}\label{2.2.3}
The inequality in Theorem \ref{tar} is sharp under natural hyperfield homomorphism $f:\mathbb{C}\rightarrow\mathbb{V}$ for polynomials of degree 2.

\begin{proof}
According to Proposition \ref{2.1.9}, here we discuss the different cases with nonzero $a$ for given polynomial $\overline{p}=T^2+c_1T+c_0$ in Poly$(\mathbb{V})$.\\
For $a\in\mathbb{V}$ with mult$_a(\overline{p})=0$, it is clear that mult$_{f^{-1}(a)}(p)=0$ for all its corresponding polynomials $p$ with respect to $f$ according to Theorem \ref{tar}. Therefore, the inequality is sharp in this case.\\
For $a\in\mathbb{V}$ with mult$_a(\overline{p})=1$, one know that $c_0=ab$ and $c_1\in a\boxplus b$ for $b\in\mathbb{V}$ and $b\neq a$, which implies that $|a-b|\leq c_1\leq a+b$. So there exists a $\theta\in[0,2\pi)$ such that $c_1=|a+be^{i\theta}|$ and then one has a corresponding polynomial
$$p(T)=T^2-(a+be^{i\theta})T+abe^{i\theta}\in\mathbb{C}[T]$$
of $\overline{p}$ with respect to $f$ to make the equality hold.\\
For $a\in\mathbb{V}$ with mult$_a(\overline{p})=2$, we know that $c_0=a^2$ and $c_1^2\leq 4c_0$ according to Proposition \ref{2.2.2}, then by calculation it is clear that $p(T)=T^2+c_1T+c_0\in\mathbb{C}[T]$ is what we want.
\end{proof}
\end{ppst}

\begin{exam}\label{2.2.4}
For any corresponding polynomial $p$ of $\overline{p}(T)=T^3+1.6T^2+0.512$ with respect to the natural hyperfield homomorphism $f:\mathbb{C}\rightarrow\mathbb{V}$, we have
\begin{center}
$\displaystyle\textup{mult}_{0.8}(\overline{p})=3>2\geq\sum_{|b|=0.8}\textup{mult}_b(p)$
\end{center}
which implies that the inequality in Theorem \ref{tar} is not sharp in this case.\\
Note that $\textup{mult}_{0.8}(\overline{p})=3$ is true since $\overline{p}\in(T+0.8)\overline{q}$ for
$$\overline{q}(T)=T^2+0.8T+0.64$$
and mult$_{0.8}(\overline{q})=2$ by Proposition \ref{2.2.2}.\\
However, one also has $\sum_{|b|=0.8}\textup{mult}_b(p)\leq2$ because for $|T|=1$, one can see
$$1.6|T|^2>|T^3|+|0.512|$$
and then by Rouche's theorem we know that $p(T)$ has exactly 2 roots inside the circle $|T|=1$ and 1 root outside of the circle $|T|=1$.
\end{exam}
\begin{rem}
Though the inequality in Theorem \ref{tar} is not always sharp, one can assert that it is sharp in some cases based on the nature of the hyperfield $\mathbb{V}$ itself with the following propositions inspired by \cite{6}.
\end{rem}

\begin{ppst}\label{2.2.6}
Let $p(T)=T^n+ c_{n-1}T^{n-1}+\cdots+ c_1T+ c_0$ be a polynomial with coefficients in $\mathbb{V}$, if the inequality
$$c_{n-k}>1+c_{n-1}+\cdots+c_{n-k+1}+c_{n-k-1}+\cdots+c_0$$
holds, then 1 is not a root of $p(T)$.

\begin{proof}
According to Proposition \ref{2.1.5}, one can suppose $c_0\neq0$ since one can always remove factors with the form $T^k$.\\
Notice that
$$c_{n-k}>1+c_{n-1}+\cdots+c_{n-k+1}+c_{n-k-1}+\cdots+c_0$$
so one can see
$$c_{n-k}\notin 1\boxplus c_{n-1}\boxplus\cdots\boxplus c_{n-k+1}\boxplus c_{n-k-1}\boxplus\cdots\boxplus c_1\boxplus c_0$$
which implies that
$$0\notin c_{n-k}\boxplus1\boxplus c_{n-1}\boxplus\cdots\boxplus c_{n-k+1}\boxplus c_{n-k-1}\boxplus\cdots\boxplus c_1\boxplus c_0$$
Since $\boxplus$ is commutative, we get $0\notin p(1)$ and then $1$ is not a root of $p$.
\end{proof}
\end{ppst}

\begin{ppst}\label{2.2.7}
Let $p(T)=T^n+ c_{n-1}T^{n-1}+\cdots+ c_1T+ c_0$ be a polynomial with coefficients in $\mathbb{V}$, if the inequality
$$c_{n-1}>1+c_{n-2}+\cdots\cdots+c_1+c_0$$
holds, then when $n\geq2$, $p$ has no roots in $\{0\}\cup[1,+\infty)$ whose multiplicity is $n$.

\begin{proof}
Suppose $p$ has a root $a$ with mult$_a(p)=n$, we want to find contradictions.\\
First, according to Proposition \ref{2.1.4}, we have $a\neq0$ otherwise $c_{i}=0$ for $i=1,2,\dots,n-1$ which leads to a contradiction. Also, from Proposition \ref{2.2.6}, this root $a$ can't be $1$.\\
When $a>1$, notice that
$$c_{n-1}a^{n-1}>(c_{n-2}+\cdots\cdots+c_1+c_0)a^{n-1}\geq c_{n-2}a^{n-2}+\cdots+ c_1a+ c_0$$
and $c_0=a^n$ (from Proposition \ref{2.1.3}), one can get
$$S=\{x:c_{n-1}a^{n-1}-(c_{n-2}a^{n-2}+\cdots+ a^n)\leq x\leq c_{n-1}a^{n-1}+c_{n-2}a^{n-2}+\cdots+a^n\}$$
where $S$ denote $c_{n-1}a^{n-1}\boxplus c_{n-2}a^{n-2}\boxplus\cdots\boxplus c_1a\boxplus a^n$. Then from
$$\begin{aligned}
c_{n-1}a^{n-1}-(a^n+c_{n-2}a^{n-2}+\cdots+ a^n)\!\!&\!&\!>&\,\,a^{n-1}(1+\sum^{n-2}_{i=0}c_i)-(a^n+\sum^{n-2}_{i=0}c_i)\\
&&=&\,\,\sum^{n-2}_{k=1}(a^{n-1}-a^i)c_i+a^{n-1}(1+a^n-2a)\\
&&>&\,\,\sum^{n-2}_{k=1}(a^{n-1}-a^i)c_i+a^{n-1}(1+a^2-2a)\\
&&=&\,\,\sum^{n-2}_{k=1}(a^{n-1}-a^i)c_i+a^{n-1}(1-a)^2\geq 0
\end{aligned}$$
one can see that $a^n\notin S$, which leads to a contradiction.
\end{proof}
\end{ppst}

\begin{ppst}\label{2.2.8}
Let $p(T)=T^n+ c_{n-1}T^{n-1}+\cdots+ c_1T+ c_0$ be a polynomial with coefficients in $\mathbb{V}$, if the inequality
$$c_{n-1}>1+c_{n-2}+\cdots\cdots+c_1+c_0$$
holds, then $p$ has no roots whose multiplicity is bigger or equal that 2 in $[1,+\infty)$.

\begin{proof}
If $p$ has no roots in $[1,\infty)$ then we are done. Now we only consider a root $a\in(1,\infty)$ of $p$ by Proposition \ref{2.2.6}. It's enough to prove that for any $q$ with $p\in(T+a)q$, $0\notin q(a)$.\\
Let $q(T)=T^{n-1}+d_{n-2}T^{n-2}+\cdots+d_0$ be such polynomial, then one can see
\begin{center}
$c_{0}=ad_{0}$ and $c_i\in ad_i\boxplus d_{i-1}$ for $i=1,\dots,n-1$
\end{center}
Therefore the inequality
$$\begin{aligned}
d_{n-2}+a&&\geq&\,\,\,c_{n-1}>1+c_{n-2}+\cdots+c_1+c_0\\
&&\geq&\quad1+|ad_{n-2}-d_{n-3}|+|ad_{n-3}-d_{n-4}|+\cdots+|ad_1-d_0|+ad_0\\
&&\geq&\quad1+ad_{n-2}-d_{n-3}+ad_{n-3}-d_{n-4}+\cdots+ad_1-d_0+ad_0\\
&&=&\quad1+ad_{n-2}+(a-1)(d_{n-2}+d_{n-3}+\dots+d_1+d_0)
\end{aligned}$$
holds, which implies that
$$1>d_{n-2}+d_{n-3}+\dots+d_1+d_0$$
Furthermore, for any $e\geq1$ we have
$$
e^{n-1}> e^{n-1}(d_{n-2}+d_{n-3}+\dots+d_1+d_0)\geq e^{n-2}d_{n-2}+e^{n-3}d_{n-3}+\dots+ed_1+d_0
$$
which leads to $e^{n-1}\notin e^{n-2}d_{n-2}\boxplus e^{n-3}d_{n-3}\boxplus\dots\boxplus ed_1\boxplus d_0$ and then $0\notin q(e)$.\\
Specifically, one can set $e=a$ to reach $0\notin q(a)$, so we are done.
\end{proof}
\end{ppst}

\begin{exam}\label{2.2.9}
Consider the polynomial
$$p(T)=T^3+5.1T+4.096$$
one has $p\in(T+1.6)q$ where $q(T)=T^2+1.6T+2.56$ since
\begin{center}
$1.6^3=4.096$, $0\in 1.6\cdot 1\boxplus 1.6$ and $5.1\in 1.6\cdot 1.6\boxplus 2.56$
\end{center}
and from Proposition \ref{2.2.2} we know that mult$_{1.6}(q)=2$ and then mult$_{1.6}(q)=3$.\\
Here one can see $c_{3-2}>1+c_2+c_0$ but mult$_{1.6}(p)=3$, which reveals that we may not be able to find a simple generalization for Proposition \ref{2.2.8}.
\end{exam}

\subsection{The Sharpness under Natural Hyperfield Homomorphism $\mathbb{C}\rightarrow\mathbb{P}$}
\begin{ppst}\label{2.3.1}
For any $x\in\mathbb{P}$, $x\boxplus x=\{x\}$ holds.

\begin{proof}
It is clear when $x=0$. If $x\neq0$, then one can see
$$x\boxplus x=\{\frac{(\alpha+\beta)x}{\|(\alpha+\beta)x\|}:\alpha,\beta\in\mathbb{R}_+\}=\{x\}$$
so we are done.
\end{proof}
\end{ppst}

\begin{exam}\label{2.3.2}
The map $f:\mathbb{C}\rightarrow\mathbb{P}$ with
\begin{center}
$re^{i\theta}\longmapsto e^{i\theta}$ when $r\neq0$ and $0\longmapsto 0$
\end{center}
is a hyperfield homomorphism.

\begin{proof}
It is clear that $f(0)=f(0),f(1)=f(1)$ and $f(xy)=f(x)f(y)$.\\
For any $xe^{i\alpha}$ and $ye^{i\beta}$, one has $f(xe^{i\alpha}\boxplus ye^{i\beta})=0$ when $x=y=0$ and
$$f(xe^{i\alpha}\boxplus ye^{i\beta})=f(\{xe^{i\alpha}+ye^{i\beta}\})=f(\{ze^{i\gamma}\})=\{e^{i\gamma}\}$$
where
$$z=\sqrt{x^2+y^2+2xy\cos(\alpha-\beta)},\,\gamma=\arctan\frac{x\sin{\alpha}+y\sin{\beta}}{x\cos{\alpha}+y\cos{\beta}}$$
when either $x$ or $y$ is nonzero. Therefore we have
\begin{center}
$f(xe^{i\alpha}\boxplus ye^{i\beta})=\{0\}\subset\{0\}=f(xe^{i\alpha})\boxplus f(ye^{i\beta})$ when $x=y=0$\\
$f(xe^{i\alpha}\boxplus ye^{i\beta})=\{\sign(x-y)e^{i\theta}\}\subset\{\pm e^{i\theta},0\}=f(xe^{i\alpha})\boxplus f(ye^{i\beta})$ when $e^{i\alpha}=-e^{i\beta}$\\
\end{center}
and
$$f(xe^{i\alpha}\boxplus ye^{i\beta})=\{e^{i\gamma}\}\subseteq\{\frac{ae^{i\alpha}+be^{i\beta}}{\|ae^{i\alpha}+be^{i\beta}\|}:a,b>0\}$$
when $xy\neq0$ since $e^{i\gamma}$ is a point on the major arc between $e^{i\alpha}$ and $e^{i\beta}$.
\end{proof}
\end{exam}

\begin{ppst}\label{2.3.3}
For $\overline{p}(T)=T^3+c_2T^2+c_1T+c_0$ in Poly$(\mathbb{P})$ and any $a\in\mathbb{P}$, there exists a corresponding polynomial $p$ with respect to the hyperfield homomorphism $f$ in Example \ref{2.3.2} such that
$$\sum_{b\in f^{-1}(a)}\textup{mult}_b(p)=\textup{mult}_a(\overline{p})$$
for $\textup{mult}_a(\overline{p})=0,1,3$.

\begin{proof}
According to Proposition \ref{2.1.8}, one only needs to prove the cases where $a\neq0$.\\
If mult$_a(\overline{p})=0$, we can let $p(T)=T^3+c_2T^2+c_1T+c_0$, then $p(a)\neq0$ follows $0\notin \overline{p}(a)$.\\
If mult$_a(\overline{p})=1$, there exists $\overline{q}(T)=T^2+d_1T+d_0$ in Poly$(\mathbb{P})$ with mult$_a(q)=0$ and
$$\overline{p}(T)\in (T-a)\boxdot (T^2+d_1T+d_0)=T^3+(-a)\boxplus d_1 T^2+(-ad_1)\boxplus d_0T-ad_0$$
which implies that $c_0=-ad_0,c_1=(-xad_1+yd_0)$ and $c_2=-za+wd_1$ for some $x,y,z,w>0$. Therefore, the following polynomial is what we want
$$\begin{aligned}
p(T)&=T^3+(-za+wd_1)T^2+\frac{wz}{x}(-xad_1+yd_0)T-\frac{wyz^2}{x}ad_0\\
&=(T-za)(T^2+wd_1T+\frac{wyz}{x}d_0)
\end{aligned}$$
If mult$_a(\overline{p})=3$, then from Theorem \ref{2.1.5} we know that $\overline{p}(T)=T^3-aT^2+a^2T-a^3$ and $p(T)=(T-a)^3$ is clearly what we want.
\end{proof}
\end{ppst}

\begin{exam} \label{2.3.4}
Consider $\overline{p}(T)=T^3+e^{i\frac{\pi}{24}}T^2+e^{i\frac{\pi}{2}}T+e^{i\frac{23\pi}{24}}$ in Poly$(\mathbb{P})$. We have
$$\overline{p}(T)\in(T+1)\boxdot(T^2+e^{i\frac{\pi}{16}}T+e^{i\frac{23\pi}{24}})$$
since $e^{i\frac{\pi}{24}}\in e^{i\frac{\pi}{2}}\boxplus e^{i\frac{\pi}{16}}$ and $e^{i\frac{\pi}{2}}\in e^{i\frac{\pi}{16}}\boxplus e^{i\frac{23\pi}{24}}$. Also, $T^2+e^{i\frac{\pi}{16}}T+e^{i\frac{23\pi}{24}}\in(T+1)(T+e^{i\frac{23\pi}{24}})$.\\
If there exists a corresponding polynomial
$$p(T)=(T+x)(T+y)(T+ze^{i\frac{23\pi}{24}})$$
of $\overline{p}$ with respect to the hyperfield homomorphism $f$ in Proposition \ref{2.3.2}, then we have
$$(x+y)+ze^{i\frac{23\pi}{24}}=me^{i\frac{\pi}{24}},\quad (x+y)ze^{i\frac{\pi}{24}}+xy=ne^{i\frac{\pi}{2}}$$
where $x,y,z,m,n>0$. One can eliminate $z$ and then has
$$(x+y)me^{i\frac{\pi}{24}}+xy-ne^{i\frac{\pi}{2}}=(x+y)^2$$
which gives
$$m(x+y)\cos{\frac{\pi}{24}}+xy-n\cos{\frac{\pi}{2}}=(x+y)^2,\quad m(x+y)\sin{\frac{\pi}{24}}=n\sin{\frac{\pi}{2}}$$
So one can see
$$x+y=\frac{n\sin{\frac{\pi}{2}}}{m\sin{\frac{\pi}{24}}}>0$$
and
$$xy=\left(\frac{n\sin{\frac{\pi}{2}}}{m\sin{\frac{\pi}{24}}}\right)^2+n\cos{\frac{\pi}{2}}-\frac{n\sin{\frac{\pi}{2}}}{\tan{\frac{\pi}{24}}}$$
Notice that $x,y$ are solutions of $t^2-(x+y)t+xy=0$, so we need $(x+y)^2\geq4xy$, which is
$$\frac{n}{m^2}\leq\frac{4}{3}\left(\frac{\sin{\frac{\pi}{24}}}{\sin{\frac{\pi}{2}}}\right)^2\left(\frac{\sin{\frac{\pi}{2}}}{\tan{\frac{\pi}{24}}}-\cos{\frac{\pi}{2}}\right)>0$$
Moreover, if $z=me^{-i\frac{22\pi}{24}}-(x+y)e^{-i\frac{23\pi}{24}}>0$, we need
$$-m\sin{i\frac{22\pi}{24}}+\frac{n\sin{\frac{\pi}{2}}}{m\sin{\frac{\pi}{24}}}\sin{i\frac{23\pi}{24}}=0\Longrightarrow \frac{n}{m^2}=\frac{\sin{\frac{22\pi}{24}}\sin{\frac{\pi}{24}}}{\sin{\frac{23\pi}{24}}\sin{\frac{\pi}{2}}}$$
Therefore, we require
$$\frac{\sin{\frac{22\pi}{24}}\sin{\frac{\pi}{24}}}{\sin{\frac{23\pi}{24}}\sin{\frac{\pi}{2}}}\leq\frac{4}{3}\left(\frac{\sin{\frac{\pi}{24}}}{\sin{\frac{\pi}{2}}}\right)^2\left(\frac{\sin{\frac{\pi}{2}}}{\tan{\frac{\pi}{24}}}-\cos{\frac{\pi}{2}}\right)$$
or equivalent
$$\frac{\sin{\frac{22\pi}{24}}}{\sin{\frac{23\pi}{24}}}\leq\frac{4}{3}\cos{\frac{\pi}{24}}$$
which is a false statement.
\end{exam}
\begin{rem}\label{2.3.5}
This example implies that the inequality in Theorem \ref{tar} is not always sharp since for this $\overline{p}$, we have
$$\textup{mult}_{-1}(\overline{p})=2>\sum_{a<0}\textup{mult}_{a}(p)$$
Furthermore, as we use the Rouche's theorem to find an example to determine the sharpness of that inequality under the natural homomorphism $\mathbb{C}\rightarrow\mathbb{V}$, it is possible for us to develop a theorem that can restrict the number of roots with certain phases for some polynomials in $\mathbb{C}[T]$ satisfying some conditions.
\end{rem}

\subsection{The Sharpness under Canonical Hyperfield Homomorphisms $\mathbb{L}\rightarrow\mathbb{K}$}
\begin{ppst}\label{2.4.1}
For every hyperfield $\mathbb{F}$, there exists a Canonical hyperfield homomorphism $f:\mathbb{F}\rightarrow\mathbb{K}$ which sends 0 to 0 and nonzero elements to nonzero elements.
\end{ppst}
\begin{rem}
It implies that $\mathbb{K}$ is an initial object in categories of hyperfields.
\end{rem}

\begin{ppst}\label{2.4.2}
For any polynomial $p$ in Poly$(\mathbb{K})$, one has
\begin{center}
$\deg p=\textup{mult}_0(p)+\textup{mult}_1(p)$.
\end{center}

\begin{proof}
From Remark \ref{2.1.7}, one can set
$$p(T)=T^n+c_{n-1}T^{n-1}+\cdots+c_1T+1$$
without losing generality. Notice that $\boxplus$ is commutative and $0\in1\boxplus1$, it is clear that
$$0\in1\boxplus1\boxplus c_{n-1}\boxplus\cdots\boxplus c_1=1\boxplus c_{n-1}\boxplus\cdots\boxplus c_1\boxplus1=p(1)$$
which implies that $p\in(T-1)q$ where $q(T)=T^{n-1}+d_{n-2}T^{n-2}+\cdots+d_1T+1$.\\
Similarly, one can get $0\in q(1)$ and, by continuing to conduct this operation inductively, we can finally get $\textup{mult}_1(p)=n$. Therefore, we are done.
\end{proof}
\end{ppst}

\begin{cor}\label{2.4.3}
For any field $\mathbb{L}$, one can find some polynomials $p$ in $\mathbb{L}[T]$ such that the equality in Theorem \ref{tar} holds with respect to the hyperfield homomorphism $f$ mentioned in Proposition \ref{2.4.1} from $\mathbb{L}$ to $\mathbb{K}$.
\end{cor}

\section{Some Existing Results and Open Questions}
\subsection{The Sharpness under Natural Hyperfield Homomorphism $\mathbb{R}\rightarrow\mathbb{S}$}
\begin{defn}\label{3.1.1}
Let $p(T)=c_n T^n+c_{n-1}T^{n-1}+\cdots+ c_1T+ c_0$ be a polynomial over $\mathbb{S}$, the \textbf{number of sign changes} in coefficients $c_n,\dots,c_0$ of $p$ is defined to be
$$\sigma(p)=\texttt{\#}\{k\geq 0|c_k=-c_{k+l+1}\neq0\textup{ and }c_{k+1}=\cdots=c_{k+l}=0\textup{ for some }l\geq0\}$$
\end{defn}
\begin{exam}
For $p(T)=T^7-T^5-T^4+T^2-1$, one has $\sigma(p)=2$.
\end{exam}
\begin{lema}\label{3.1.3}Let $p(T)=c_n T^n+ c_{n-1}T^{n-1}+\cdots+c_1T+ c_0$ be a polynomial over $\mathbb{S}$, then one has mult$_1(p)=\sigma(p)$.
\end{lema}
\begin{rem}\label{3.1.4}
This lemma cited from \cite{3} gives us a good way to compute the multiplicity of roots in hyperfield of signs. In fact, one can also determine the multiplicity of $-1$ and $0$ by it with simple deduction.
\end{rem}
\begin{ppst}\label{3.1.5}
For every sequence $(s_i)_{i=0}^n\subset\mathbb{S}$ with $s_n>0$ and $s_0\neq0$ , there exist a corresponding polynomial $p(T)=c_nT^n+\cdots+c_1T+c_0\in\mathbb{R}[T]$ such that
\begin{center}
$\displaystyle\sum_{a>0}\textup{mult}_a(p)=\sigma(p)$
\end{center}
and $\sign(c_i)=s_i$ for $i=0,\dots,n$.

\begin{proof}
The proof is given by David J. Grabiner in \cite{5} where he gives a polynomial
$$p(T)=\sum_{i=0}^n\frac{s_i}{(2n)^{i^2}}T^i$$
which satisfies that $\sum_{a>0}\textup{mult}_a(p)=\sigma(p)$. It is valid because, for two every consecutive $s_k,s_l$ with $s_ks_l<0$ in such polynomial, one can check that $p((2n)^{2k})p((2n)^{2l})<0$ and then know that $p$ has a root between $(2n)^{2k}$ and $(2n)^{2l}$.
\end{proof}
\end{ppst}
\begin{rem} \label{3.1.6}
Here one knows that there always exist a $p(T)\in\mathbb{R}[T]$, such that
$$\sum_{a\in f^{-1}(1)}\textup{mult}_a(p)=\sigma(p)=\sigma(\overline{p})=\textup{mult}_1(\overline{p})$$
for any $\overline{p}(T)$ with coefficients in $\mathbb{S}$. Notice that $f(a)=0$ is equal to $a=0$, so one has
$$\textup{mult}_0(p)=\textup{mult}_0(\overline{p})$$
Additionally, for every consecutive $s_k,s_l$ with $s_ks_l>0$, one has $p(-(2n)^{2k})p(-(2n)^{2l})<0$ and then get a negative root, so
$$\sum_{a\in f^{-1}(-1)}\textup{mult}_a(p)=\textup{mult}_{-1}(\overline{p})$$
where $\textup{mult}_{-1}(\overline{p})=\textup{mult}_{1}(\overline{q})$ where $q(T):=p(-T)$ in Proposition \ref{3.1.5}.\\
Therefore, the inequality in Theorem \ref{tar} under $f:\mathbb{R}\rightarrow\mathbb{S}$ is always sharp.
\end{rem}
\subsection{The Sharpness under Valuative Hyperfield Homomorphisms $\mathbb{R}\rightarrow\mathbb{T}$}
\begin{defn}\label{3.2.1}
Given a field $\mathbb{K}$, a map $\nu:\mathbb{K}\rightarrow\mathbb{R}\cup\{-\infty\}$ is called a \textbf{(Krull) valuation} if\\
(1) $\nu^{-1}(-\infty)=0$;\\
(2) $\nu(ab)=\nu(a)+\nu(b)$\\
(3) $\nu(a+b)\leq\max\{\nu(a),\nu(b)\}$
for all $a,b\in\mathbb{K}$.
\end{defn}
\begin{exam}
Consider the map
$$f:\mathbb{R}\rightarrow\mathbb{R}\cup\{0\},x\longmapsto \log{|x|}$$
where $f(0)=\log{0}=-\infty$. Then it is clear that $f$ is a valuation.
\end{exam}
\begin{ppst}\label{3.2.3}
A map $\nu:\mathbb{K}\rightarrow\mathbb{R}\cup\{-\infty\}$ is a valuation if and only if it is a homomorphism of hyperfields.
\end{ppst}
\begin{thm}\label{3.2.4} (Fundamental theorem for the tropical hyperfield) Given a polynomial $p=T^n+c_{n-1}T^{n-1}+\cdots+c_1T+c_0$ in Poly$(\mathbb{T})$, then\\
(1) There is a unique sequence $a_1,\dots,a_n\in\mathbb{T}$, up to permutation of the indices, such that
$$c_{n-i}\in\bigboxplus_{e_1<\cdots<e_i}\!\!\!\!_{\mathbb{T}}\,\,a_{e_1}\cdots a_{e_{i}}$$
for all $i=1,\dots,n$ where $\{e_1,\dots,e_i\}$ is a subset of $\{1,2,\dots,n\}$ and $\boxplus_{\mathbb{T}}$ is the hyperaddition in $\mathbb{T}$.\\
(2) The equalities
$$\textup{mult}_a(p)=\#\{i\in\{1,\dots,n\}\,|\,a=a_i\}$$
hold for all $a\in\mathbb{T}$.
\end{thm}
\begin{cor}
Given a polynomial $p$ in $\textup{Poly}(\mathbb{T})$, one has $\deg(p)=\sum_{a\in\mathbb{T}}\textup{mult}_a(p)$.
\end{cor}
\begin{rem}
This theorem cited from \cite{3} is critical to the study of tropical hyperfield. In fact, the techniques used to prove this theorem can be also applied to prove similar results about general valuative hyperfields. Moreover, one can immediately get the sharpness of the inequality in Theorem \ref{tar} under all (valuative) hyperfield homomorphisms $\tau:\mathbb{C}\rightarrow\mathbb{F}$ where $\mathbb{F}$ is a valuative hyperfield.
\end{rem}

\newpage
\subsection{Some Open Questions}
In this part, we list some open questions associated with this paper.\\

\noindent\textbf{Question 1}\label{q1} As is discussed in Remark \ref{2.3.5}, we wonder if we can develop a theorem that can restrict the number of roots with certain phases for some polynomials in $\mathbb{C}[T]$ satisfying some conditions.

\begin{defn}\label{3.3.1}
A hyperfield $\mathbb{F}$ is called stringent if $a\neq-b$ implies $a\boxplus b$ is a singleton for all $a,b\in\mathbb{F}$.
\end{defn}
\begin{exam}
It is clear that $\mathbb{S}$ and all valuative hyperfields are stringent.
\end{exam}
\noindent\textbf{Question 2}\label{q2}
In Nathan Bowler and Rudi Pendavingh's joint work \cite{4}, stringent hyperfields can be classified into three types which are associated with a field $\mathbb{F}$, the Krasner hyperfield $\mathbb{K}$ and the hyperfield of signs $\mathbb{S}$ respectively. Since the inequality is sharp in both the three hyperfields, we wonder if it is true that the inequality is sharp in all stringent hyperfields.\\

\noindent\textbf{Question 3}\label{q3}
In \cite{3}, Baker and Oliver describe Poly$(\mathbb{F})$ in terms of ordered blueprints and they give the matroids associated with ordered blueprints in other works of them. Therefore, we are curious about whether we can connect matroid theory with the inequality in Theorem \ref{tar}.\\

\noindent\textbf{Question 4}\label{q4}
Can we give an algebraic proof for the sharpness under $\mathbb{R}\rightarrow\mathbb{S}$?

\end{document}